\documentclass[12pt]{amsart}

\usepackage{a4,amsmath,amsfonts,amssymb,amsthm,bbold,array,graphicx,overpic}

\DeclareMathOperator{\convOp}{conv} \DeclareMathOperator{\SigOp}{S}
\DeclareMathOperator{\weightOp}{w} 
 \DeclareMathOperator{\AOp}{A}
\DeclareMathOperator{\BOp}{B} \DeclareMathOperator{\littleOOp}{o}
\DeclareMathOperator{\bigOOp}{O}

\newcommand{\zerovec}{\ensuremath{\mathbb{0}}}
\newcommand{\onevec}{\ensuremath{\mathbb{1}}}
\newcommand{\uvec}[1]{\ensuremath{\mathbb{e}_{#1}}}
\newcommand{\setdef}[2]{\{{#1}\,:\,{#2}\}}
\newcommand{\setdefbig}[2]{\big\{{#1}\,:\,{#2}\big\}}
\newcommand{\setdefadap}[2]{\left\{{#1}\,:\,{#2}\right\}}
\newcommand{\lerlex}{\prec_{\text{\tiny rlex}}}
\newcommand{\leqrlex}{\preceq_{\text{\tiny rlex}}}
\newcommand{\symdiff}{\bigtriangleup}
\newcommand{\zo}[1]{\ensuremath{\{0,1\}^{#1}}}
\newcommand{\zostar}[1]{\ensuremath{\{0,1\}^{#1}\setminus\{\zerovec\}}}

\newcommand{\compX}[1]{\ensuremath{X^{\prec {#1}}}}
\newcommand{\compP}[1]{\ensuremath{P^{\prec {#1}}}}
\newcommand{\compXblock}[2]{\ensuremath{\compX{#1}_{#2}}}
\newcommand{\compPblock}[2]{\ensuremath{\compP{#1}_{#2}}}

\newcommand{\knapX}[1]{\ensuremath{X^{< {#1}}}}
\newcommand{\knapP}[1]{\ensuremath{P^{< {#1}}}}

\newcommand{\knapPblock}[2]{\ensuremath{\knapP{#1}_{#2}}}

\newcommand{\sig}[2]{\ensuremath{\sigma_{#2}({#1})}}
\newcommand{\weight}[1]{\ensuremath{\weightOp({#1})}}
\newcommand{\Sig}[1]{\ensuremath{\SigOp({#1})}}
\newcommand{\SigGr}[2]{\ensuremath{\SigOp^{>{#2}}({#1})}}
\newcommand{\coSigGr}[2]{\ensuremath{\overline{\SigOp}^{>{#2}}({#1})}}

\newcommand{\coSig}[1]{\ensuremath{\overline{\SigOp}({#1})}}
\newcommand{\Apq}[3]{\ensuremath{\AOp^{\prec{#1}}_{#2,#3}}}
\newcommand{\Apqx}[4]{\ensuremath{\Apq{#1}{#2}{#3}({#4})}}
\newcommand{\Bpq}[3]{\ensuremath{\BOp^{\prec{#1}}_{#2,#3}}}
\newcommand{\Bpqx}[4]{\ensuremath{\Bpq{#1}{#2}{#3}({#4})}}
\newcommand{\R}{\ensuremath{\mathbb{R}}}
\newcommand{\Q}{\ensuremath{\mathbb{Q}}}
\newcommand{\N}{\ensuremath{\mathbb{N}}}
\newcommand{\rangez}[1]{\ensuremath{[{#1}]_0}}
\newcommand{\range}[1]{\ensuremath{[{#1}]}}
\newcommand{\bfz}{\ensuremath{\mathbf{0}}}
\newcommand{\bfo}{\ensuremath{\mathbf{1}}}
\newcommand{\bfst}{\ensuremath{\mathbf{\star}}}

\newcommand{\lb}[2]{g_{\text{nfac}}({#1},{#2})}
\newcommand{\lbe}[2]{g_{\text{avdeg}}({#1},{#2})}
\newcommand{\littleO}[1]{\littleOOp({#1})}
\newcommand{\bigO}[1]{\bigOOp({#1})}

\newcommand {\defi} {:=}  

\newcommand{\Qplus}{\Q^{\geq 0}}
\newcommand{\net}[1]{\mathcal{N}\!\left({#1}\right)}
\newcommand{\expans}[1]{\mathcal{X}\!\left({#1}\right)}

\newtheorem{definition}{Definition}
\newtheorem{proposition}{Proposition} 
\newtheorem{theorem}{Theorem}

\begin{document}

\title{ Revlex-Initial 0/1-Polytopes} \thanks{This work has been
  supported by the DFG Research Group \emph{Algorithms, Structure,
    Randomness} and by the DFG Research Center \textsc{Matheon}}
\author{Volker Kaibel} \author{Rafael Mechtel}
\email{kaibel@zib.de} \email{mechtel@math.tu-berlin.de}
\address{Zuse Institute Berlin\\Takustr.~7\\14195~Berlin\\ Germany}
\address{TU Berlin\\ MA 6--2\\ Stra\ss e des 17. Juni 136\\ 10623
  Berlin \\Germany}

\date{July 22, 2005} 

\maketitle

\begin{abstract}
  We introduce revlex-initial 0/1-polytopes as the convex hulls of
  reverse-lexicographically initial subsets of 0/1-vectors. These
  polytopes are special knapsack-polytopes. It turns out that they
  have remarkable extremal properties. In particular, we use these
  polytopes in order to prove that the minimum numbers $\lb{d}{n}$ of
  facets and the minimum average degree $\lbe{d}{n}$ of the graph of a
  $d$-dimensional 0/1-polytope with~$n$ vertices satisfy $\lb{d}{n}\le
  3d$ and $\lbe{d}{n}\le d+4$. We furthermore show that, despite the
  sparsity of their graphs, revlex-initial 0/1-polytopes satisfy a
  conjecture due to Mihail and Vazirani, claiming that the graphs of
  0/1-polytopes have edge-expansion at least one.
\end{abstract}


\section{Introduction}

Let us call a subset~$X$ of $\{0,1\}^d$ \emph{revlex-initial} if, for
every $x\in X$, all points in $\{0,1\}^d$ that are
reverse-lexicographically smaller than~$x$ are contained in~$X$. The
convex hulls of revlex-initial subsets of $\{0,1\}^d$ are the
\emph{revlex-initial 0/1-polytopes}. Phrased differently, the
revlex-initial 0/1-polytopes are the convex hulls of those sets of
0/1-vectors of length~$d$ that correspond to the binary
representations of all numbers $0,1,\dots,n-1$ for some~$n$. In
particular, for every $1\le n\le 2^d$ there is precisely one
revlex-initial 0/1-polytope with~$n$ vertices in~$\R^d$.

Why should one be interested in such special polytopes? The general
interest in 0/1-polytopes stems from their importance in combinatorial
optimization. Investigations of 0/1-polytopes like traveling salesman
polytopes, cut polytopes, stable set polytopes, and matching polytopes
have not only led to beautiful insights into the interplay of
combinatorics and geometry, but also to great algorithmic progress
with respect to the corresponding optimization problems. From that
work on such special 0/1-polytopes quite a few general questions on
0/1-polytopes have emerged, such as, e.g., the question for the
maximal number of facets a $d$-dimensional 0/1-polytope may have (see
Ziegler~\cite{Zie00}).

With respect to this extremal question, B\'ar\'any and
P\'or~\cite{BP01} obtained a remarkable result. They showed that a
random $d$-dimensional 0/1-polytope with roughly $2^{{d}/{\log_2 d}}$
vertices in expectation has at least (roughly) $2^{({1}/{4})d\log_2
  d}$ facets. Recently, this bound was even improved to
$2^{({1}/{2})d\log_2 d}$ by Gatzouras, Giannopoulos, and
Markoulakis~\cite{GGM04}. The best known upper bound currently is
$\bigO{(d-2)!}$ (due to Fleiner, Kaibel, and Rote~\cite{FKR00}).  It
turns out that the revlex-initial 0/1-polytopes studied in this paper
give some answers to two reverse extremal questions: How \emph{few}
facets or edges can a $d$-dimensional 0/1-polytope with a specified
number of vertices have?

Note that, somewhat different to the class of general polytopes, the
number of vertices of a 0/1-polytope may impose severe restrictions on
the combinatorial type. For instance, a 0/1-polytope is simple if and
only if it is the product of (0/1-)simplices (Kaibel and
Wolff~\cite{KW00}). Thus, $d$-dimensional simple 0/1-polytopes
with~$n$ vertices do only exist if there is a factorization $n=\prod
n_i$ of~$n$ with $d=\sum(n_i-1)$ .  Therefore, within the realm of
0/1-polytopes, it seems interesting to investigate extremal questions
for \emph{all} (reasonable) pairs $(d,n)$.

Our paper contains three main results.

(1) Revlex-initial 0/1-polytopes in $\R^d$ have no more than $3d$
facets (Theorem~\ref{thm:facets}); from this we deduce that the
smallest number of facets $\lb{d}{n}$ of a $d$-dimensional
0/1-polytope with exactly~$n$ vertices satisfies $\lb{d}{n}\le 3d$ for
all~$d$ and~$n$ and $\lb{d}{n(d)}\le d+\littleO{d}$ if $n(d)$ grows
sub-exponentially with~$d$ (Theorem~\ref{thm:lbfac}).

(2) The average degree of every revlex-initial 0/1-polytope in $\R^d$
is at most $d+4$ (Theorem~\ref{thm:edgenumcompr}); from this we deduce
that the smallest average degree $\lbe{d}{n}$ of a $d$-dimensional
0/1-polytope with exactly~$n$ vertices satisfies $\lbe{d}{n}\le d+4$
(Theorem~\ref{thm:edgenumbers}).

Since revlex-initial 0/1-polytopes have extremely sparse graphs, at
first sight they look like candidates for counter-examples to an
important conjecture due to Mihail and Vazirani (cited, e.g.,
in~\cite{FM92,Mih92}) stating that the graph of every 0/1-polytope has
edge-expansion at least one.  However, supporting that conjecture, we
prove:

(3) Revlex-initial 0/1-polytopes have edge-expansion at least one
(Theorem~\ref{thm:edge-expansion}); from this we deduce that, for
every (reasonable) pair $(d,n)$, there are $d$-dimensional
0/1-polytopes with~$n$ vertices, sparse graphs, and edge-expansion at
least one (Theorem~\ref{thm:expall}).

The context in which we came to study the special class of
revlex-initial 0/1-polytopes is described in
Section~\ref{subsec:incr}.  They appeared from investigating an
apparently strange behavior of certain convex hull algorithms on
random 0/1-polytopes.

The notion of revlex-initial subsets of $\{0,1\}^d$, or, equivalently,
of a system of subsets of $\{1,\dots,d\}$, is not new.  It is related
to the notion of \emph{compression} of a set system, which plays an
important role in the Kruskal-Katona theorem (see, e.g.,
\cite[Thm.~8.32]{Zie98}) characterizing the $f$-vectors of simplicial
complexes. Here, a system $\mathcal{S}$ of subsets of $\{1,\dots,d\}$
(corresponding to a subset $X\subseteq\{0,1\}^d$) is called
\emph{compressed} if, for every~$i$, the subsystem of~$\mathcal{S}$
containing all sets from~$\mathcal{S}$ of cardinality~$i$ is
reverse-lexicographically initial within the $i$-subsets of
$\{1,\dots,d\}$.  Clearly, every revlex-initial subset of $\{0,1\}^d$
corresponds to a compressed system of subsets of $\{1,\dots,d\}$, but
the converse is not true.

In the context of the Kruskal-Katona theorem only compressed set
systems that are closed under taking subsets are considered. Of
course, all revlex-initial 0/1-polytopes correspond to compressed set
systems with that property (i.e., revlex-initial 0/1-polytopes are
\emph{monotone}).  But even more: Exploiting the interpretation in
terms of binary representations of numbers, one finds that
revlex-initial 0/1-polytopes are a special kind of knapsack polytopes
(see Section~\ref{sec:def}).

Note that the terminus 'compressed polytope' has already been coined
with a different meaning (see, e.g., \cite{Sta80}).

\subsection*{Acknowledgments} We are thankful to Jens Hillmann for
computer implementations and for performing several computer
experiments and to Michael Joswig for stimulating
discussions. Furthermore we wish to thank two anonymous referees for
their helpful remarks as well as Marc E. Pfetsch and G\"unter
M. Ziegler for carefully reading an earlier version of the
manuscript. 


\section{Definitions}
\label{sec:def}

Throughout the paper, we assume that $d$ is a positive integer number.
We start with fixing some notions and notation. 

\begin{definition}[Index ranges]
  For a positive integer number~$k$, let
  $$
  \range{k}\ :=\ \{1,2,\dots,k\}\quad\text{and}\quad\rangez{k}\ :=\ 
  \{0,1,\dots,k-1\}\ .
  $$
  We will identify $\R^d$ with $\R^{\rangez{d}}$, i.e. vectors
  $x\in\R^d$ have components $x_0$, $x_1$, \dots, $x_{d-1}$, similarly
  for $\N^d$. 
\end{definition}

\begin{definition}[Reverse-lexicographical order]
  A point $x\in\{0,1\}^d$ is \emph{reverse-lexicographically smaller}
  than another point $y\in\{0,1\}^d\setminus\{x\}$ ($x\lerlex y$) if
  $x_{i_{\max}}<y_{i_{\max}}$ holds for
  $i_{\max}:=\max\setdef{i}{x_i\not= y_i}$.  We denote $x\leqrlex y$
  if $x=y$ or $x\lerlex y$ hold for $x,y\in\{0,1\}^d$.
\end{definition}

For $x\in\{0,1\}^d$ denote $\Sig{x}:=\setdef{i\in\rangez{d}}{x_i=1}$.
Then we have
$$
x\lerlex y\ \Leftrightarrow\ \max(\Sig{x}\symdiff
\Sig{y})\,\in\,\Sig{y}
$$
for all $x,y\in\{0,1\}^d$ ($x\not= y$), where $\symdiff$ denotes
the symmetric difference of two sets.

\begin{definition}[Revlex-Initial 0/1-polytope]
  A subset $X\subseteq\zo{d}$ is \emph{ revlex-initial} if, for every
  $x\in X$, it contains all $y\in\zo{d}$ with $y\lerlex x$.
  For $v \in \{0,1\}^d$ define
  $$
  \compX{v}\ \defi\ \setdef{x\in\zo{d}}{x\lerlex v}\ .
  $$
  A 
  \emph{ revlex-initial 0/1-polytope} is the convex
  hull of any 
   revlex-initial 0/1-set. We denote
  $$
  \compP{v}\ \defi\ \convOp \compX{v}\ .
  $$
\end{definition}

Since $\lerlex$ defines a total ordering of $\zo{d}$, every 
 revlex-initial 0/1-set $X$ with $\lvert X \rvert < 2^d$ is of
the form $\compX{v}$ for some $v\in\zostar{d}$. 
Note that $v\not\in\compP{v}$.

\begin{definition}[Signature of a 0/1-point]
  Let $v\in\zostar{d}$. Its \emph{weight} $\weight{v}\defi\onevec^\top
  v$ is the number of ones of~$v$. Its \emph{signature} is the vector
  $$
  (\sig{v}{1},\dots,\sig{v}{\weight{v}})
  $$
  with
  $$
  \Sig{v}\ =\ \{\sig{v}{1},\dots,\sig{v}{\weight{v}}\}
  \quad\text{and}\quad
  \sig{v}{1}>\sig{v}{2}>\dots>\sig{v}{\weight{v}}\ .
  $$
  Further we define the index set of all zero-components  
  $$
  \coSig{v}\defi\rangez{d}\setminus\Sig{v}.
  $$
\end{definition}

\begin{definition}[Block decomposition]
  For a 0/1-point $v\in\zostar{d}$ with signature
  $(\sig{v}{1},\dots,\sig{v}{\weight{v}})$, we call
  $$
  \compXblock{v}{q}\ \defi\ 
  \setdef{x\in\zo{d}}{ %
    x_{\sig{v}{q}}=0, x_{\sig{v}{q}+1}=v_{\sig{v}{q}+1}, \ldots, %
    x_{d-1} = v_{d-1} }
  $$
  (for $q\in\range{\weight{v}}$) the \emph{blocks} of~$\compP{v}$.
  Clearly, $\compX{v}$ is the disjoint union
  $$
  \compX{v}=\compXblock{v}{1}\uplus\dots\uplus\compXblock{v}{\weight{v}}
  $$
  of its blocks. The faces
  $\compPblock{v}{q}\defi\convOp\compXblock{v}{q}$ are the \emph{block
    faces} of $\compP{v}$. The vector
  $$
  (\dim\compPblock{v}{1},\dots,\dim\compPblock{v}{\weight{v}}) \ =\ 
  (\sig{v}{1},\dots,\sig{v}{\weight{v}})
  $$
  is the \emph{signature} of the 
   revlex-initial 0/1-polytope
  $\compP{v}$.
\end{definition}

\begin{table}[htdp]
\caption{Example illustrating some of the definitions: We have $d=10$,
  $\weight{v}=5$, $\Sig{v}=\{0,2,3,6,9\}$, and
  $\coSig{v}=\{1,4,5,7,8\}$.}  
\centerline{%
$
\renewcommand{\arraystretch}{1.5}
\newlength{\tmplength}
\settowidth{\tmplength}{\sig{v}{4}}
\begin{tabular}{l>{\centering}p{\tmplength}>{\centering}p{\tmplength}>{\centering}p{\tmplength}>{\centering}p{\tmplength}>{\centering}p{\tmplength}>{\centering}p{\tmplength}>{\centering}p{\tmplength}>{\centering}p{\tmplength}>{\centering}p{\tmplength}c}
{v} &   \bfo & \bfz & \bfo & \bfo & \bfz & \bfz & \bfo  & \bfz & \bfz & \bfo \\
\hline
\text{indices} &  0 & 1 & 2 & 3 & 4 & 5 & 6 & 7 & 8 & 9 \\
\text{signature}  &  \sig{v}{5} & & \sig{v}{4} & \sig{v}{3} & & & \sig{v}{2} & & & \sig{v}{1} \\
\hline
\compPblock{v}{1} & \bfst &\bfst &\bfst &\bfst &\bfst &\bfst &\bfst &\bfst &\bfst &\bfz \\
{\compPblock{v}{2}} & \bfst &\bfst &\bfst &\bfst &\bfst &\bfst &\bfz &\bfz &\bfz &\bfo \\
{\compPblock{v}{3}} & \bfst &\bfst &\bfst &\bfz &\bfz &\bfz &\bfo &\bfz &\bfz &\bfo \\
{\compPblock{v}{4}} & \bfst &\bfst &\bfz &\bfo &\bfz &\bfz &\bfo &\bfz &\bfz &\bfo \\
{\compPblock{v}{5}} & \bfz &\bfz &\bfo &\bfo &\bfz &\bfz &\bfo &\bfz &\bfz &\bfo 
\end{tabular}
$
}
\end{table}%

As mentioned in the introduction, \emph{} revlex-initial
0/1-polytopes are a special kind of knapsack polytopes. Indeed, for $d
\in \N$ we define $a \in \N^d$ as $a_i \defi 2^i$.  Then for two
0/1-vectors $v,w \in \zo{d}$ we have $v \lerlex w$ if and only if
$a^\top v < a^\top w$ holds.  
Thus we can identify each natural number $n \in \N$ with a unique 0/1-vector
$v \in \zo{d}$ for a unique $d$ such that $n = a^\top v$ and $v_{d-1} = 1$.
Therefore we write $\knapP{n}$ with $n \in \N$ instead of $\compP{v}$
with $v \in \zo{d}$ with $v_{d-1} = 1$.
With the above identification, $\knapP{n}$ has exactly the $n$
vertices corresponding to the numbers $0, 1, \ldots, n-1$.
In other words, $\compP{v}$ with $v \in \zo{d}$ is the knapsack
polytope $\convOp \setdef{x \in \zo{d}}{a^\top x \le a^\top v -1}$.

\section{The Facets of  Revlex-Initial 0/1-Polytopes}

\subsection{Optimizing Linear Functions}

For $c\in\R^d$ and $I\subseteq\rangez{d}$, define
$$
c^+(I)\ \defi\ \sum_{i\in I}\max\{c_i,0\}\ .
$$
The following statement follows immediately from the block
decomposition of revlex-initial 0/1-polytopes.

\begin{proposition}
\label{prop:opt}
For every $v\in\zo{d}\setminus\{\zerovec\}$ and $c\in\R^d$, we have
$$
\max\setdef{c^\top x}{x\in\compP{v}} \ =\ 
\max\setdefbig{\sum_{p=1}^{q-1}c_{\sig{v}{p}}+c^+(\rangez{\sig{v}{q}})}%
{q \in \range{\weight{v}}}\ .
$$
\end{proposition}

In particular, the optimization problem
$\max\setdef{c^\top x}{x\in\compP{v}}$ (for given $v\in\zo{d}$ and
$c\in\Q^d$) can be solved in polynomial time.

\subsection{A Linear Description}

If $i\in\coSig{v}$ and
$x\in\compX{v}$ with $x_i=1$, then $x_j=0$ must hold for some
$j\in\Sig{v}$ with $j>i$. Let us denote
$$
\SigGr{v}{i}\ \defi\ \setdef{j\in\Sig{v}}{j>i}\ \text{ and }
\coSigGr{v}{i}\ \defi\ \setdef{j\in\coSig{v}}{j>i}\ .
$$
We will use similar notations with respect to~$<$, $\le$ and~$\ge$.
Thus, the inequalities
\begin{equation}
\label{eq:rlexieq}
x_i+\sum_{j\in\SigGr{v}{i}}x_j\ \le\ |\SigGr{v}{i}|\qquad
\text{for all }i\in\coSig{v}
\end{equation}
and (since $v\not\in\compP{v}$)
\begin{equation}
\label{eq:rlexieqall}
\sum_{j\in\Sig{v}}x_j\ =\ v^\top x\ =\ \le\ |\Sig{v}|-1
\end{equation}
are valid for~$\compP{v}$.
These inequalities are minimal cover inequalities. In fact, they are
all minimal cover inequalities of the knapsack polytope $\compP{v}$.

\begin{theorem}[Linear descriptions of  revlex-initial 0/1-polytopes]
  \label{thm:lindescr}
  For every $v\in\zostar{d}$ the 
  revlex-initial 0/1-polytope~$\compP{v}$ has the following linear
  description: 
  \begin{equation} \label{eq:lindescr}
    \compP{v}\ =\ \setdef{x\in\R^d}{%
      \zerovec\le x\le\onevec\,,
      x\text{ satisfies~\eqref{eq:rlexieq} and~\eqref{eq:rlexieqall}}}
  \end{equation}
\end{theorem}

\begin{proof}
  Denote the polytope defined by the right-hand side
  of~\eqref{eq:lindescr} by $Q(v)$. Thus, $Q(v)$ is the set of all
  $x\in\R^d$ satisfying the following system of inequalities:
\begin{eqnarray}
\label{eq:ieq1}
-x_i & \le &  0  \qquad\text{for all }i\in\rangez{d}\\
\label{eq:ieq2}
 x_i & \le &  1 \qquad\text{for all }i\in\rangez{d}\\
\label{eq:ieq3}
x_i+\sum_{j\in\SigGr{v}{i}}x_j & \le &  |\SigGr{v}{i}| \qquad\text{for
  all }i\in\coSig{v} \\
\label{eq:ieq4}
\sum_{j\in\Sig{v}}x_j & \le & \weight{v}-1
\end{eqnarray}
Denote by $A$ the matrix with the left-hand side coefficients of the
inequalities in \eqref{eq:ieq3} and \eqref{eq:ieq4}. The rows of $A$
can be put into an order, such that $A$ is an interval matrix.
Thus $A$ is total unimodular (see, e.g., \cite[Example~7,
p.~279]{Sch86}), and appending the identity matrices $I_d$ and $-I_d$
does not change total unimodularity.
                                
Since the righthand sides of the inequalities in
\eqref{eq:ieq1}--\eqref{eq:ieq4} are integers, all vertices of $Q(v)$
are integer vectors and by the inequalities of type \eqref{eq:ieq1}
and \eqref{eq:ieq2} they are binary vectors.
Therefore $Q(v) = \compP{v}$, since $Q(v) \cap \zo{d} = \compX{v}$.
\end{proof}

\subsection{The Facet Defining Inequalities}

Let us first describe the dimension of a revlex-initial 0/1-polytope.

\begin{proposition} \label{prop:dim}
  For each $v\in\zostar{d}$ the dimension of the 
   revlex-initial 0/1-polytope~$\compP{v}$ is
  $$
  \dim\compP{v}\ =\
  1+\max\left(\setdef{i\in\rangez{d}}{\uvec{i}\lerlex v} \cup \{-1\}\right)\ ,
  $$
  
  In our knapsack notation we have for $n \in \N$
  $$
  \dim\knapP{n}\ =\ 1+\max\setdef{i\in\N \cup \{-1\}}{2^i<n} =
  \min\setdef{j\in\N}{n\le2^j}\ .
  $$
\end{proposition}

\begin{proof}
  This follows from the block decomposition of~$\compP{v}$.
\end{proof}

In particular, $\compP{v}$ is full-dimensional if and only if
$\uvec{d-1}\lerlex v$ (that is, $2^{d-1} < n \le 2^d$). The following
three propositions describe the facets of full-dimensional
revlex-initial 0/1-polytopes.

\begin{proposition}
\label{prop:facetszero}
For each $v\in\zo{d}$ with $\uvec{d-1}\lerlex v$ and for every
$i\in\rangez{d}$, the inequality $x_i\ge 0$ defines a facet
of~$\compP{v}$.
\end{proposition}

\begin{proof}
  By Theorem~\ref{thm:lindescr}, the inequalities
  \eqref{eq:ieq1}--\eqref{eq:ieq4} provide a linear description
  of~$\compP{v}$. Since the trivial inequalities~\eqref{eq:ieq1} are
  the only ones in this description which have negative coefficients,
  none of them can be conically combined from others. Hence, they all
  define facets of~\compP{v} (since~\compP{v} is full-dimensional).
\end{proof}

\begin{proposition}
\label{prop:Sigfacet}
For each $v\in\zo{d}$ with $\uvec{d-1}\lerlex v$, the inequality
$\sum_{j\in\Sig{v}}x_j\le \weight{v}-1$ defines a facet
of~$\compP{v}$.
\end{proposition}

\begin{proof}
  The inequality $\sum_{j\in\Sig{v}}x_j\le \weight{v}-1$ is the only
  inequality in the linear description
  \eqref{eq:ieq1}--\eqref{eq:ieq4} of~$\compP{v}$ provided by
  Theorem~\ref{thm:lindescr} that is violated by the point $v$, which
  is not contained in~$\compP{v}$. Thus, that inequality must define a
  facet of~$\compP{v}$.
\end{proof}

\begin{proposition}
\label{prop:facetsone}
For each $v\in\zo{d}$ with $\uvec{d-1}\lerlex v$ and for every
$i\in\rangez{d}$, the inequality $x_i\le 1$ defines a facet
of~$\compP{v}$ unless
\begin{center}
  $\weight{v}=2$ and $i\in\Sig{v}$\\
  or\\
  $\sig{v}{2}<d-2$ and $\sig{v}{2}<i\le d-1$
\end{center} 
(in which cases they do not define facets).
\end{proposition}

\begin{proof}
  Unless one of the exceptions listed in the proposition holds, all
  inequalities from the linear description
  \eqref{eq:ieq1}--\eqref{eq:ieq4} of~$\compP{v}$ provided by
  Theorem~\ref{thm:lindescr} that have a positive $i$-th coefficient
  have right-hand-side at least two. Since the only ones with negative
  $i$-th coefficient have right-hand-side zero, the inequality $x_i\le
  1$ cannot be conically combined from the others in that linear
  description. Hence it defines a facet of~\compP{v} (since~\compP{v}
  is full-dimensional).
  
  In case of $\weight{v}=2$ and $i\in\Sig{v}$, let~$j$ be such that
  $\Sig{v}=\{i,j\}$. Thus, $x_i\le 1$ is the sum of
  inequality~\eqref{eq:ieq4} and $-x_j\le 0$. Hence, it does not
  define a facet of~\compP{v}.
  
  Finally, consider the case $\sig{v}{2}<d-2$. If $\sig{v}{2}<i<d-1$,
  then the type-\eqref{eq:ieq3} inequality $x_i+x_{d-1}\le 1$ implies
  $x_i\le 1$ by adding $-x_{d-1}\le 0$. If $i=d-1$ then the
  type-\eqref{eq:ieq3} inequality $x_j+x_{d-1}\le 1$ for any
  $\sig{v}{2}<j<d-1$ implies $x_{d-1}\le 1$ by adding $-x_j\le 0$.
  Thus, in both cases, $x_i\le 1$ does not define a facet
  of~\compP{v}.
\end{proof}

\begin{proposition}
\label{prop:coSigfacets}
For each $v\in\zo{d}$ with $\uvec{d-1}\lerlex v$ and for every
$i\in\coSig{v}$, the inequality $x_i+\sum_{j\in\SigGr{v}{i}}x_j\le
|\SigGr{v}{i}|$ defines a facet of~$\compP{v}$ unless
$i<\sig{v}{\weight{v}}$ (in which case it does not define a facet).
\end{proposition}

\begin{proof}
  For each $i\in\coSig{v}$ with $i>\sig{v}{\weight{v}}$, the
  inequality $x_i+\sum_{j\in\SigGr{v}{i}}x_j\le |\SigGr{v}{i}|$ is the
  only inequality in the linear description
  \eqref{eq:ieq1}--\eqref{eq:ieq4} of~$\compP{v}$ 
  that is violated by the point
  $v+\uvec{i}-\uvec{\sig{v}{\weight{v}}}$, which is not contained
  in~$\compP{v}$. Thus, that inequality must define a facet
  of~$\compP{v}$.
  
  If $i<\sig{v}{\weight{v}}$, then $x_i+\sum_{j\in\SigGr{v}{i}}x_j\le
  |\SigGr{v}{i}|$ does not define a facet since it equals the sum of
  the two inequalities $\sum_{j\in\Sig{v}}x_j\le \weight{v}-1$ and
  $x_i\le 1$.
\end{proof}

Combining Theorem~\ref{thm:lindescr} and the five preceding
propositions, we obtain the following result.

\begin{theorem}[Facets of  revlex-initial 0/1-polytopes]
\label{thm:facets}
Let $v\in\zo{d}$ with $\uvec{d-1}\lerlex v$, i.e., \compP{v} is a
full-dimensional revlex-initial 0/1-polytope. Let
$$
D(v)\ \defi\ D_1(v)\cup D_2(v)
$$
with
$$
D_1(v)\defi\begin{cases}
  \Sig{v} & \text{if }\weight{v}=2\\
  \varnothing & \text{otherwise}
\end{cases}
$$
and
$$
D_2(v)\defi\begin{cases}
  \{\sig{v}{2}+1,\dots,d-1\} & \text{if } \sig{v}{2}<d-2\\
  \varnothing & \text{otherwise}
          \end{cases}\ .
          $$
\begin{enumerate}
\item The following system is a minimal (with respect to. inclusion) linear
  description of~\compP{v} by facet defining inequalities:
  $$
\begin{array}{rcll}
x_i & \ge &  0  & \text{for all }i\in\rangez{d}\\
x_i & \le &  1   & \text{for all }i\in\rangez{d}\setminus D(v)\\
x_i+\displaystyle\sum_{j\in\SigGr{v}{i}}x_j & \le &  |\SigGr{v}{i}| &
\text{for all }i\in\coSig{v}, i>\sig{v}{\weight{v}}\\ 
\displaystyle\sum_{j\in\Sig{v}}x_j & \le & \weight{v}-1
\end{array} \ .
$$
\item The number of facets of~\compP{v} is
  $$
  f_{d-1}(\compP{v})\ =\ 
  2d+\big|\setdef{\sig{v}{\weight{v}}<i<\sig{v}{2}}{v_i=0}\big|
  +\epsilon\ ,
  $$
  where
  $$
  \epsilon\ \defi\ 
\begin{cases}
  -1   & \text{if }\weight{v}=2\\
  0  & \text{if }\weight{v}>2, v_{d-2}=0\\
  1 & \text{otherwise (i.e., $\weight{v}>2, v_{d-2}=1$)}
\end{cases}\ .
$$
We have
$$
2d-1\le f_{d-1}(\compP{v})\le 3d-2\ .
$$
The minimum number $2d-1$ of facets is attained if and only if
$\weight{v}=2$, and the maximum $f_{d-1}(\compP{v})=3d-2$ is achieved
only by $v=\uvec{0}+\uvec{d-2}+\uvec{d-1}$ (for $d\ge 3$).
\end{enumerate}

\end{theorem}

\begin{figure}[ht]
  \centering
  \includegraphics[width=.6\textwidth]{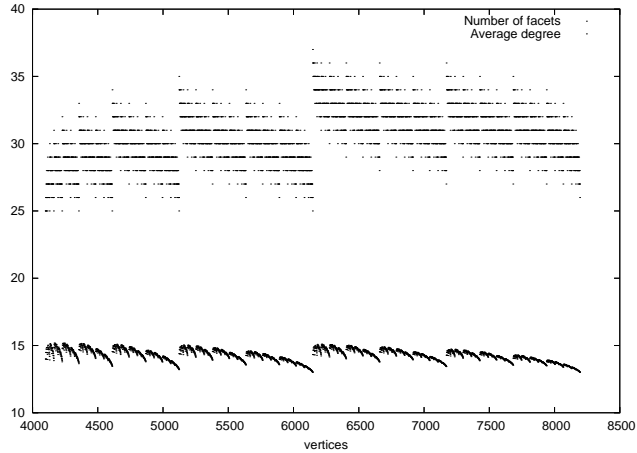}
  \caption{
    The numbers of facets and the average degrees of all
    full-dimensional revlex-initial 0/1-polytopes for $d=13$.
  }
  \label{fig:nfacavdeg-13}
\end{figure}

See Figure~\ref{fig:nfacavdeg-13} for an illustration of the facet numbers
of  revlex-initial 0/1-polytopes.

\subsection{Incremental Convex-Hull Algorithms}
\label{subsec:incr}

\begin{figure}[ht]
  \centering
  \includegraphics[width=.6\textwidth]{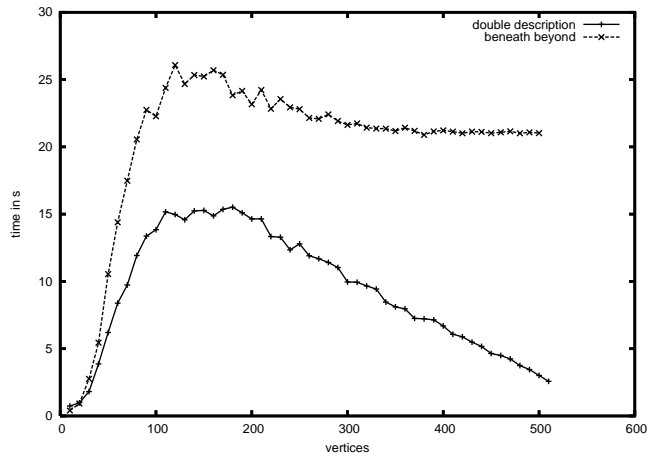}
  \caption{
    Incremental convex hull algorithms: running times on
    $9$-dimensional random 0/1-polytopes.
  }
  \label{fig:incremental}
\end{figure}

The origin of our investigations on revlex-initial 0/1-polytopes lies
in some experiments on computing the convex hulls of random
0/1-polytopes that we performed with the \texttt{polymake} system.
Some of the results of the experiments are illustrated in
Figure~\ref{fig:incremental}, showing the running times for computing
the convex hulls of (uniformly) random 0/1-polytopes in~$\R^d$
depending on the number~$n$ of vertices. The picture shows two curves,
one for the \emph{beneath-beyond} and one for the
\emph{double-description} method (where \texttt{polymake} uses
Komei Fukuda's implementation \texttt{cdd} for the latter method).

These two methods are \emph{incremental} in the sense that they
iteratively compute the convex hull of the first $i+1$ vertices from
the convex hull of the first~$i$ vertices. Since $n-1$ vertices of a
random 0/1-polytope with~$n$ vertices should make a random
0/1-polytope with $n-1$ vertices, we had expected the curves to be
monotonically increasing. However, the \emph{first} $n-1$ vertices do
only make a (uniform) random 0/1-polytope with $n-1$ vertices if the
order of the~$n$ vertices is (uniformly) random.

As it turned out, this is not the case for random 0/1-polytopes
produced by the \texttt{polymake} system. Instead, the \texttt{rand01}
client of \texttt{polymake} is implemented in such a way that the
vertices of the random 0/1-polytope produced appear in lexicographic
order. This led us to studying revlex-initial 0/1-polytopes.

And in fact, our results on the facet numbers of  revlex-initial
0/1-polytopes make the curves in Figure~\ref{fig:incremental}
plausible: For 0/1-polytopes with large numbers of vertices, which
furthermore are lexicographically ordered, the intermediate polytopes
appearing during the runs of incremental convex hull algorithms are
quite close to  revlex-initial 0/1-polytopes. Therefore, it is
plausible that these intermediate polytopes have extremely few facets 
compared to random 0/1-polytopes with the same numbers of vertices.  

In particular, if the $2^d$ vertices of the entire cube are ordered
lexicographically then the total number of facets of all intermediate
polytopes produced by an incremental convex hull algorithm to compute
the cube is bounded from above by $3d\cdot 2^d$, while for an
arbitrary (even for a random) ordering there might be intermediate
polytopes with super-exponentially many vertices (due to the results
of B\'ar\'any and P\'or~\cite{BP01} and Gatzouras, Giannopoulos, and
Markoulakis~\cite{GGM04}).

These results indicate that it might be a good strategy to sort the
vertices lexicographically before applying an incremental convex hull
algorithm to a 0/1-polytope. However, we do not yet have any thorough
computational study to support this.


\section{The Graphs of  Revlex-Initial 0/1-Polytopes}

\subsection{Characterization of Adjacency}

The one-dimensional faces of a polytope (forming its
\emph{$1$-skeleton} or \emph{graph}) are particularly important, for
instance, since the simplex algorithm for linear programming proceeds
along them. Moreover, in the special case of 0/1-polytopes, the graphs
are important also for different reasons (see
Section~\ref{subsec:edgeexp}).

Here, we describe the graphs of revlex-initial 0/1-polytopes.

\begin{definition}
  For $v\in\zostar{d}$ and $1\le p<q\le\weight{v}$ and $x\in
  \zo{d}$ we define the sets
  \begin{align*}
    \Apqx{v}{p}{q}{x}\ \defi\ \big\{ z\in\compPblock{v}{p} \;:\;&
      z_i=x_i\text{ for all }0\le i<\sig{v}{q},\\
      & z_{\sig{v}{q}}=0,\,z_{\sig{v}{r}}=1\text{ for all }p<r<q
    \big\}
  \end{align*}
  and
  \begin{align*}
    \Bpqx{v}{p}{q}{x}\ \defi\ \big\{z\in\compPblock{v}{p} \;:\;z &
      z_i=x_i\text{ for all }0\le i<\sig{v}{q}, \\
      & z_{\sig{v}{q}}=1,\,z_{\sig{v}{r}}=1\text{ for all }p<r<q 
    \big\}\ .
  \end{align*}
\end{definition}

\begin{table}[htdp]
\caption{Illustration of the definitions (with $p=1$, $q=4$, and
  $x=(1,0,x_2,x_3,x_4,x_5,x_6,x_7,x_8,x_9) \in
  \zostar{10}$) on the example from Section~\ref{sec:def}.} 
\centerline{%
$
\renewcommand{\arraystretch}{1.5}
\settowidth{\tmplength}{\sig{v}{4}}
\begin{tabular}{l>{\centering}p{\tmplength}>{\centering}p{\tmplength}>{\centering}p{\tmplength}>{\centering}p{\tmplength}>{\centering}p{\tmplength}>{\centering}p{\tmplength}>{\centering}p{\tmplength}>{\centering}p{\tmplength}>{\centering}p{\tmplength}c}
{v} &   \bfo & \bfz & \bfo & \bfo & \bfz & \bfz & \bfo  & \bfz & \bfz & \bfo \\
\hline
\text{indices} &  0 & 1 & 2 & 3 & 4 & 5 & 6 & 7 & 8 & 9 \\
\text{signature}  & \sig{v}{5}  & & \sig{v}{4} & \sig{v}{3} & & & \sig{v}{2} & & & \sig{v}{1} \\
\hline
\compPblock{v}{1} & \bfst &\bfst &\bfst &\bfst &\bfst &\bfst &\bfst &\bfst &\bfst &\bfz \\
\Apqx{v}{1}{4}{x}     & \bfo &\bfz &\bfz &\bfo &\bfst &\bfst &\bfo &\bfst &\bfst &\bfz \\
\Bpqx{v}{1}{4}{x}     & \bfo &\bfz &\bfo &\bfo &\bfst &\bfst &\bfo &\bfst &\bfst &\bfz \\
\compPblock{v}{4}  & \bfst &\bfst &\bfz &\bfo &\bfz &\bfz &\bfo &\bfz &\bfz &\bfo 
\end{tabular}
$
}
\end{table}%

\begin{theorem}[Graphs of  revlex-initial 0/1-polytopes]
  \label{thm:graph}
  For $v\in\zostar{d}$, the graph of the corresponding revlex-initial
  0/1-polytope $\compP{v}$ has the following structure.
  \begin{enumerate}
  \item Let $x\in\compX{v}$ be a vertex of~$\compP{v}$ contained in the
    block~$\compPblock{v}{q}$. 
    Let $p$ be some block number with $1\le p<q$.
    \begin{enumerate}
    \item The vertex $x$ is adjacent to all vertices of
      $\Apqx{v}{p}{q}{x}$.
    \item If
      $\max(\setdef{i\in\rangez{\sig{v}{q}}}{x_i\not=v_i}\cup\{-1\})\not\in\Sig{v}$
      then~$x$ is also adjacent to all vertices of
      $\Bpqx{v}{p}{q}{x}$.
    \end{enumerate}
  \item The graph of $\compP{v}$ does not contain any other edges than
    the (cube-)edges of the blocks
    $\compPblock{v}{1}$,\dots,$\compPblock{v}{\weight{v}}$ and the ones
    described in part~(i) of this theorem.
  \end{enumerate}
\end{theorem}

\begin{proof}
  For the proof of part~(1), let us denote by~$F$ the face
  of~$\compP{v}$ that is defined by the following equations:
\begin{eqnarray}
\label{eq:F1}
z_i=x_i               & \qquad(0\le i<\sig{v}{q}) \\
\label{eq:F2}
z_{\sig{v}{r}}=1 & \qquad(p<r<q)\\
\label{eq:F3}
z_i=v_i              & \qquad(\sig{v}{p}<i)\
\end{eqnarray}

The claim in~(a) follows from the fact that the only vertices of the
face $\setdef{z\in F}{z_{\sig{v}{q}}=0}$ of~$\compP{v}$ are the
vertices of $\Apqx{v}{p}{q}{x}$ and~$x$ itself. Since
$\Apqx{v}{p}{q}{x}$ is contained in the hyperplane defined by
$z_{\sig{v}{p}}=0$, while~$x$ is not, that face must be the 
pyramid with base $\Apqx{v}{p}{q}{x}$ and apex~$x$.

In order to prove part~(b), assume
$\max(\setdef{i\in\rangez{\sig{v}{q}}}{x_i\not=v_i}\cup\{-1\})\not\in\Sig{v}$.
Thus, there is no block $\compPblock{v}{r}$ with $r>q$ that has a
common vertex with the face $F$.  Hence, the only vertices of that
face are the vertices of $\Apqx{v}{p}{q}{x}$, $\Bpqx{v}{p}{q}{x}$,
and~$x$ itself. Again, since $\Apqx{v}{p}{q}{x}$ and
$\Bpqx{v}{p}{q}{x}$ are contained in the hyperplane defined by
$z_{\sig{v}{p}}=0$, while~$x$ is not, that face must be the pyramid
with base $\convOp(\Apqx{v}{p}{q}{x}\cup\Bpqx{v}{p}{q}{x})$ and
apex~$x$.

For the proof of part~(2), suppose that~$x$ and~$y$ are adjacent
vertices of~$\compP{v}$ not contained in the same block. We may assume
$x\in\compPblock{v}{q}$ and $y\in\compPblock{v}{p}$ with $1\le p<
q\le\weight{v}$.
 
We will first show that~$y$ is contained in the face~$F$
of~$\compP{v}$ defined in the proof of part~(1). Therefore, we have to
prove that \eqref{eq:F1}--\eqref{eq:F3} is satisfied by $z=y$.
 
Let us assume~\eqref{eq:F1} is not satisfied by $z=y$, i.e., there is
some $0\le i<\sig{v}{q}$ with $x_i\not=y_i$. If we denote, for
$a,b\in\zo{d}$, by $a\oplus b$ the component-wise addition modulo two,
then we have $x\oplus\uvec{i}\in\compP{v}$ (since $i<\sig{v}{q}$) and
$y\oplus\uvec{i}\in\compP{v}$ (since $i<\sig{v}{q}<\sig{v}{p}$) with
$$
\{x\oplus\uvec{i},y\oplus\uvec{i}\}\ \ne\ \{x,y\}
$$
(since $x_{\sig{v}{p}}=1\ne0=y_{\sig{v}{p}}$). But then
$$
\tfrac{1}{2}(x+y)\ =\ \tfrac{1}{2}(x\oplus\uvec{i}+y\oplus\uvec{i})
$$
contradicts the adjacency of~$x$ and~$y$. Thus, $z=y$
satisfies~\eqref{eq:F1}.
 
If~\eqref{eq:F2} would not be satisfied by $z=y$, then there was some
$p<r<q$ with $y_{\sig{v}{r}}=0$. Due to $x\in\compPblock{v}{q}$,
$x_{\sig{v}{r}}=1$ holds. Thus, we have
$x-\uvec{\sig{v}{r}}\in\compP{v}$ and
$y+\uvec{\sig{v}{r}}\in\compP{v}$ (since $y\in\compPblock{v}{p}$ with
$r<p$). Again,
$$
\{x-\uvec{\sig{v}{r}},y+\uvec{\sig{v}{r}}\}\ \ne\ \{x,y\}
$$
holds, and therefore,
$$
\tfrac{1}{2}(x+y)\ =\ 
\tfrac{1}{2}\big((x-\uvec{\sig{v}{r}})+(y+\uvec{\sig{v}{r}})\big)
$$
contradicts the adjacency of~$x$ and~$y$. Hence, \eqref{eq:F2} is
satisfied by $z=y$.
 
Since $q>p$ and $x\in\compPblock{v}{q}$, $y\in\compPblock{v}{p}$, we
clearly have $x_i=y_i=v_i$ for all $i>\sig{v}{p}$. Therefore,
also~\eqref{eq:F3} is satisfied by $z=y$, and thus, the claim $y\in F$
is proved.
 
We obtain
$y\in\Apqx{v}{p}{q}{x}\cup\Bpqx{v}{p}{q}{x}$. It hence suffices to
show that, in case of $y\in\Bpqx{v}{p}{q}{x}$, we have
$$
\max(\setdef{i\in\rangez{\sig{v}{q}}}{x_i\ne v_i}\cup\{-1\})\ 
\not\in\ \Sig{v}\ .
$$
Therefore, suppose we have $y\in\Bpqx{v}{p}{q}{x}$ and there is
some $q<s\le\weight{v}$ with $x_{\sig{v}{s}}=0$ and $x_i=v_i$ for all
$\sig{v}{s}<i<\sig{v}{q}$.  Then we have
$y-\uvec{\sig{v}{q}}\in\compP{v}$ (due to $y\in\compPblock{v}{p}$,
$p<q$, and $y_{\sig{v}{q}}=1$) and $x+\uvec{\sig{v}{q}}\in\compP{v}$
(in fact: $x+\uvec{\sig{v}{q}}\in\compPblock{v}{s}$). Also here, we
have $$
\{ x+\uvec{\sig{v}{q}},y-\uvec{\sig{v}{q}}\}\ \ne\ \{x,y\}\ ,
$$
and thus,
$$
\tfrac{1}{2}(x+y)\ =\ 
\tfrac{1}{2}\big((x+\uvec{\sig{v}{q}})+(y-\uvec{\sig{v}{q}})\big)
$$
contradicts the adjacency of~$x$ and~$y$.
 
\end{proof}

\subsection{The Number of Edges}

Having the structural description given in Theorem~\ref{thm:graph} at
hand, we can now derive a formula for the number of edges of a
revlex-initial 0/1-polytope.

\begin{theorem}[Edge numbers of  revlex-initial 0/1-polytopes]
\label{thm:edgenumcompr}
  For $v\in\zostar{d}$, the graph of the corresponding 
   revlex-initial 0/1-polytope $\compP{v}$ has
   $$
   \sum_{p=1}^{\weight{v}}2^{\sig{v}{p}}\Big(\tfrac{\sig{v}{p}}{2}
   + \sum_{q=p+1}^{\weight{v}}2^{p-q}\Big(2-
   \Big(\sum_{r=q+1}^{\weight{v}}2^{\sig{v}{r}}\Big)
   2^{-\sig{v}{q}}\Big)\Big) 
  $$
  edges. In particular, its average node degree is bounded from
  above by $d+4$.
\end{theorem}

\begin{proof}
  The statement on the average degree follows from the exact
  expression for the number of edges: Inside the (outermost) brackets,
  the fraction $\tfrac{\sig{v}{p}}{2}$ is bounded from above by
  $\tfrac{d}{2}$ while the remaining sum clearly is at most~$2$. Thus
  the number of edges is at most $(\tfrac{d}{2}+2)$ times the number
  $\sum 2^{\sig{v}{p}}$ of vertices of~$\compP{v}$.
 
  In order to determine the total number of edges, let $1\le
  p<q\le\weight{v}$. We have
  $$
  \dim\Apqx{v}{p}{q}{x}\ =\ \dim\Bpqx{v}{p}{q}{x}\ =\ 
  (p+\sig{v}{p})-(q+\sig{v}{q})\ =:\ \delta_{p,q}
  $$
  for each $x\in\compXblock{v}{q}$.
  
  Clearly, the number of edges between~$\compPblock{v}{q}$
  and~$\compPblock{v}{p}$ described in part~(1a) of
  Theorem~\ref{thm:graph} thus is
  $$
  2^{\sig{v}{q}}\cdot 2^{\delta_{p,q}}\ =\ 2^{p+\sig{v}{p}-q}\ .
  $$
  The number of $x\in\compXblock{v}{q}$ that do not satisfy the
  condition of part~(1b) of Theorem~\ref{thm:graph} is
  $\sum_{r=q+1}^{\weight{v}}2^{\sig{v}{r}}$.  Thus, the number of
  edges between~$\compPblock{v}{q}$ and~$\compPblock{v}{p}$ described
  in part~(1b) is
  $$
  \Big(2^{\sig{v}{q}} - \sum_{r=q+1}^{\weight{v}}2^{\sig{v}{r}}
  \Big) \cdot 2^{\delta_{p,q}}
  \ =\ 2^{p+\sig{v}{p}-q}-\Big(\sum_{r=q+1}^{\weight{v}}2^{\sig{v}{r}}\Big) 2^{\delta_{p,q}}\ .
  $$
  
  Therefore, the total number of edges is
  $$
  \sum_{p=1}^{\weight{v}}\sig{v}{p}2^{\sig{v}{p}-1}\ +\ \sum_{1\le
    p<q\le\weight{v}}
  \bigg(2\cdot2^{p+\sig{v}{p}-q} - 
   \Big(\sum_{r=q+1}^{\weight{v}}2^{\sig{v}{r}}\Big)
  2^{(p+\sig{v}{p})-(q+\sig{v}{q})}\bigg)\ ,
  $$
  where the first sum accounts for the edges inside the blocks and
  the second one (the double-sum) counts the edges running across
  different blocks. That expression can easily be simplified to the
  one stated in the theorem.
 
\end{proof}

\subsection{The Edge-Expansion}
\label{subsec:edgeexp}

The geometry of a 0/1-polytope~$P$ (more precisely: its $1$-skeleton,
i.e., its graph) defines a natural neighborhood structure on the set
system~$\mathcal{S}$ corresponding to the vertices of~$P$. Such a
neighborhood structure can be used in order to design random walk
algorithms for generating elements from~$\mathcal{S}$ at random
(according to a certain pre-specified probability distribution).
Random walk algorithms are of great importance, for instance with
respect to randomized approximative counting algorithms (see, e.g.,
\cite{JS97}).

In many cases, the neighborhood structure defined geometrically via
the associated 0/1-polytope has turned out to be quite appropriate for
designing such random walk algorithms. A crucial parameter with
respect to the time complexity of these methods is the \emph{edge
  expansion} of the neighborhood structure. The rule of thumb here is
that the expansion should be bounded from below polynomially in
$1/d$ (where $d$ is the dimension of the polytope) in order
to achieve an efficient time algorithm. 

\begin{definition}
  The edge expansion $\expans{G}$ of a graph $G = (V,E)$ is defined as
  \begin{align}
    \expans{G} \defi 
    & \min\setdefadap{\frac{|\delta(S)|}{|S|}}{S\subset V,\
      0<|S|\leq\frac{|V|}{2}}\nonumber 
  \end{align}
  (with~$\delta(S)$ denoting the set of all edges with one end node
  in~$S$ and the other one in~$V\setminus S$).
\end{definition}

It has been conjectured by Mihail and Vazirani (cited, e.g.,
in~\cite{FM92,Mih92}) that the graph of every 0/1-polytope has edge
expansion at least one.  In fact, this conjecture is known to be true
for several classes of 0/1-polytopes, including stable set polytopes,
(perfect) matching polytopes, and polytopes associated with the bases
of \emph{balanced} (in particular: regular) matroids
(see~\cite{Mih92,FM92,Kai04}). For more details and references, we
refer to~\cite{Kai04}.

Here, further supporting the Mihail-Vazirani conjecture, we prove that
despite the sparsity of their graphs, revlex-initial 0/1-polytopes
have edge expansion at least one.

\begin{theorem}[Edge expansion of  revlex-initial 0/1-polytopes] 
  \label{thm:edge-expansion}
  For $v\in\zostar{d}$, the graph of the corresponding revlex-initial
  0/1-polytope $\compP{v}$ has edge expansion at least one.
\end{theorem}

In order to bound the edge expansion of a graph~$G=(V,E)$ from below
we will construct certain flows in the (uncapacitated)
network~$\net{G}=(V,A)$, where~$A$ contains for each edge~$\{u,v\}\in
E$ both arcs~$(u,v)$ and~$(v,u)$.  This strategy dates back to the
method of ``canonical paths'' developed by Sinclair
(see~\cite{Sin93}). The extension to flows was explicitly exploited by
Morris and Sinclair~\cite{SM99}. Feder and Mihail~\cite{FM92} use
random canonical paths, which can equivalently be formulated in terms
of flows.

The crucial idea is to construct for each ordered pair $(x,y)\in
V\times V$ a flow $\phi_{(x,y)}:A\longrightarrow\Qplus$ in the
network~$\net{G}$ sending one unit of some commodity from~$x$ to~$y$.
Define the multi-commodity flow (MCF) $\phi\defi\sum_{(x,y)\in V\times 
  V}\phi_{(x,y)}$ as the sum of all the flows $\phi_{(x,y)}$. 
By
$$
\phi_{\max}\defi\max\setdef{\phi(a)}{a\in A}
$$
we denote the maximal amount of $\phi$-flow on any arc. By
construction of~$\phi$, the total amount $\phi(S:V{\setminus} S)$ of
$\phi$-flow leaving~$S$ is at least $|S|\cdot(n-|S|)$, where $n=|V|$.
On the other hand, we have $\phi(S:V{\setminus} S)\leq
\phi_{\max}\cdot|\delta(S)|$. This implies
$|S|\cdot(n-|S|)\leq\phi_{\max}\cdot|\delta(S)|$, and hence, if
$|S|\leq\frac{n}{2}$ holds,
$$
\frac{|\delta(S)|}{|S|}\geq \frac{n}{2\cdot\phi_{\max}}\enspace.
$$
Thus, we have proven
\begin{equation}
  \label{eq:expphimax}
  \expans{G}\geq\frac{n}{2\cdot\phi_{\max}}\enspace.
\end{equation}

In the light of inequality~(\ref{eq:expphimax}) it is clear that the
task is to construct a flow~$\phi$ as above with $\phi_{\max} \le
\tfrac n2$ (where $n = \lvert V \rvert$).

\begin{proof}[Proof of Theorem \ref{thm:edge-expansion}]
  We will use the notations $\knapP{n}\defi\compP{v}\subset\R^{d}$ and
  $\knapX{n}\defi\compX{v}$, where~$n\in\N$ is the number having
  binary representation~$v$ (i.e., $n$ is the number of vertices
  of~$\knapP{n}=\compP{v}$). Clearly, we may assume $v_{d-1}=1$, i.e.,
  $n>2^{d-1}$ and $\dim\knapP{n}=d$. Thus, in particular, the
  dimension~$d$ and the 0/1-vector $v\in\zo{d}$ are uniquely
  determined by the vertex number~$n$.
                                
  We will prove the theorem by showing via induction on~$n$ that, for
  every $n \in \N$, there is an MCF $\phi^n = \sum_{(x,y) \in
    \knapX{n} \times \knapX{n}} \phi^n_{(x,y)}$ on
  $\net{G(\knapP{n})}$ such that $\phi^n_{\max} \le \tfrac n2$.
  
  The statement obviously holds for $n = 2$, since in that case, the
  polytope $\knapP{n}$ consists of two vertices joint by an edge.
  
  Thus let us suppose that for all $2\le n' < n$ there is such an MCF
  $\phi^{n'}$ on $\net{G(\knapP{n'})}$ with $\phi^{n'}_{\max} \le
  \tfrac {n'}{2}$. The induction step, i.e., the construction of an
  appropriate MCF $\phi^n$, will be subdivided into two cases.
  
  Let $G \defi G(\knapP{n})$. For a subset $A$ of the nodes of~$G$, we
  denote by $G[A]$ the subgraph of~$G$ induced by~$A$ (similarly, we
  use $\net{G}[A]$). Two 0/1-polytopes~$P$ and~$Q$ are called
  \emph{0/1-equivalent} if they can be transformed into each other by
  (potentially) lifting one of them into the space of the other and
  applying a symmetry of the cube (i.e., by flipping and permuting
  coordinates). Of course, such a transformation induces an
  isomorphism between the graphs of~$P$ and~$Q$.  Note that for $w \in
  \zo{d}$ the vertex $w \oplus \uvec{0}$ is the
  one obtained by flipping the first coordinate of $w$.

  \paragraph{Case 1 ($v_0 = 0$)}
  Define the following faces of $\knapP{n}$ and the corresponding vertex sets:
  \begin{align*}
    F_A \defi \setdef{w \in \knapP{n}}{w_0 = 0} &&
    F_B \defi \setdef{w \in \knapP{n}}{w_0 = 1} \\
    X_A \defi \setdef{w \in \knapX{n}}{w_0 = 0} &&
    X_B \defi \setdef{w \in \knapX{n}}{w_0 = 1}\ .
  \end{align*}
  Then, for every $x \in X_A$, we have $x \oplus \uvec{0} \in X_B$
  (and vice versa). Thus $\knapP{n}$ is a prism over $F_A$. In
  particular, $F_A$ and $F_B$ are 0/1-equivalent. Furthermore, they
  both are 0/1-equivalent to $\knapP{n'}$ with $n' = \tfrac n2$. Thus,
  $G[X_A]$ and $G[X_B]$ both are isomorphic to $G(\knapP{n'})$.
                                
  Let $\phi^A$ and $\phi^B$ be the MCFs induced by $\phi^{n'}$ on
  $\net{G}[X_A]$ and $\net{G}[X_B]$, respectively. Thus $\phi^A_{\max}
  = \phi^B_{\max} = \phi^{n'}_{\max} \le \tfrac{n}{4}$ by the
  induction hypothesis.
  Now we construct the MCF $\phi^n$ on $\net{G}$ by defining each
  $\phi^n_{(x,y)}$ in the following way (note that $G[X_A]$ and
  $G[X_B]$ are edge-disjoint):
  $$  
  \begin{array}{ll}
    x,y \in X_A : & \phi^n_{(x,y)} \defi \phi^A_{(x,y)} \\
    x,y \in X_B : & \phi^n_{(x,y)} \defi \phi^B_{(x,y)} \\
    x \in X_A, y \in X_B : &
    \phi^n_{(x,y)} \defi \Psi_{(x,x \oplus \uvec{0})} + \phi^B_{(x \oplus
      \uvec{0},y)}\\ 
    x \in X_B, y \in X_A : & 
    \phi^n_{(x,y)} \defi \Psi_{(x,x \oplus \uvec{0})} + \phi^A_{(x \oplus
      \uvec{0},y)}. 
  \end{array}
  $$
  Here, $\Psi_{(x,x \oplus \uvec{0})}$ denotes the flow just
  sending one unit along the arc $(x,x \oplus \uvec{0}) \in A$ (and
  nothing along any other arc).
  In the resulting MCF $\phi^n$, every arc $(x,x \oplus \uvec{0})$
  with $x \in X_A$ carries one unit of flow for each of the $\lvert
  X_B \rvert = \tfrac n2$ pairs $(x,y)$, $y \in X_B$.  The same holds
  for the reverse arcs $(x,x \oplus \uvec{0})$, $x \in X_B$.  Thus we
  have $\phi(x,x \oplus \uvec{0}) = \tfrac n2$ for every such arc, and
  we conclude
  $$
  \phi^n_{\max} \le \max\left\{ \tfrac n2, 2 \cdot \phi^{n'}_{\max} 
  \right\} = \frac n2.
  $$

  \paragraph{Case 2 ($v_0 = 1$)} Let $\hat{x}\in\zo{d}$ be the 
  revlex-predecessor of~$v$, i.e., the 0/1-vector corresponding to the
  number $(n-1)$. Then $\hat{x}$ is the ``last'' vertex of $\knapP{n}$
  and $\{\hat{x}\} = \knapPblock{n}{\weight{n}}$ is the ``last'' block
  of $\knapP{n}$.  Define the following faces of $\knapP{n}$ and the
  corresponding vertex sets (see Fig.~\ref{fig:case2}):
  \begin{align*}
    F_A &\defi \setdef{w \in \knapP{n}}{w_0 = 0} &
    F_B &\defi \setdef{w \in \knapP{n}}{w_0 = 1} \\
    X_A &\defi \setdef{w \in \knapX{n}}{w_0 = 0} &
    X_B &\defi \setdef{w \in \knapX{n}}{w_0 = 1} \cup \{\hat{x}\} \\
    X'_A &\defi X_A \setminus \{\hat{x}\} &
    X'_B &\defi X_B \setminus \{\hat{x}\}\ .
  \end{align*}
  Thus $\knapP{n}$ is a \emph{partial prism} over $F_A$, i.e.,
  $\knapP{n}$ arises from a true prism over $F_A$ via removing the
  vertex $(\hat{x}+\uvec{0})$ corresponding to $\hat{x}$ (and taking
  the convex hull).
  
  \begin{figure}[ht]
    \centering
    \begin{overpic}[scale=.5,tics=5]{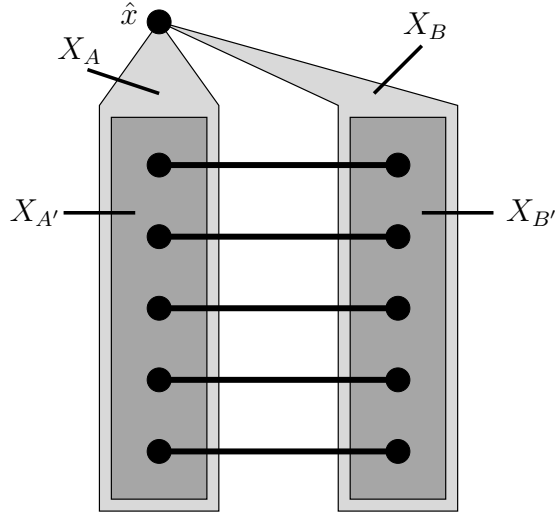}
      \put(12,97){$\hat{x}$}
      \put(-.5,90){$X_A$}
      \put(-10,58){$X_{A'}$}
      \put(68,95){$X_B$}
      \put(88.5,58){$X_{B'}$}
    \end{overpic}
    \caption{
      Illustration of the sets used in case~2 of the proof of
      Theorem~\ref{thm:edge-expansion}.
    }
    \label{fig:case2}
  \end{figure}
  
  We will first prove that there is a spanning subgraph of $G[X_B]$
  that is isomorphic to $G[X_A]$. Indeed, this is a simple consequence
  of the fact that $\knapP{n}$ is a partial prism over $F_A$ (with
  $\hat{x}$ being the ``not duplicated vertex''): Every edge
  $\{a,a'\}$ of $F_A$ with $a,a'\not=\hat{x}$ gives rise to a
  quadrangular $2$-face $\{a,a',a\oplus\uvec{0},a'\oplus\uvec{0}\}$ of
  the partial prism (showing that
  $\{a\oplus\uvec{0},a'\oplus\uvec{0}\}$ is an edge of~$G$), and every
  edge $\{\hat{x},a\}$ of $F_A$ yields a triangular $2$-face
  $\{\hat{x},a,a\oplus\uvec{0}\}$ of the partial prism (showing that
  $\{\hat{x},a\oplus\uvec{0}\}$ is an edge of~$G$).
  
  Hence, there is a spanning subgraph of $G[X_B]$ that is isomorphic
  to $G[X_A]$. Furthermore, the face $F_A$ of~$\knapP{n}$ is
  0/1-equivalent to $\knapP{n'}$ with $n' = (n+1)/2$.  Therefore, the
  MCF $\phi^{n'}$ induces MCFs $\phi^A$ and $\phi^B$ on $\net{G}[X_A]$
  and $\net{G}[X_B]$, respectively, with
  $$
  \phi^A_{\max} = \phi^B_{\max} = \phi^{n'}_{\max} \le \tfrac{n+1}{4}
  $$ 
  by the induction hypothesis.  
  
  With $\alpha \defi \tfrac {n-1}{n+1} < 1$ we have
  $(1+\alpha)\phi^{n'}_{\max} \le \tfrac n2$. Thus we can increase
  each of the flows $\phi^A$ and $\phi^B$ by an $\alpha$-fraction
  without making the flow exceed the desired limit
  of~$n/2$ at any arc.
  We construct the MCF $\phi^n$ on $\net{G}$ by defining each
  $\phi^n_{(x,y)}$ in the following way (note that $G[X_A]$ and
  $G[X_B]$ are edge-disjoint):
  $$
  \begin{array}{ll}
    x,y \in X_A : & \phi^n_{(x,y)} \defi \phi^A_{(x,y)} \\
    x,y \in X_B : & \phi^n_{(x,y)} \defi \phi^B_{(x,y)} \\
    x \in X'_A, y \in X'_B : &
    \phi^n_{(x,y)} \defi \alpha\left( \Psi_{(x,x \oplus \uvec{0})} +
      \phi^B_{(x \oplus \uvec{0},y)} \right) 
    + (1-\alpha) \left( \phi^A_{(x,\hat{x})} + \phi^B_{(\hat{x},y)}
    \right) \\
    x \in X'_B, y \in X'_A  : & 
    \phi^n_{(x,y)} \defi \alpha\left( \Psi_{(x,x \oplus \uvec{0})} +
      \phi^A_{(x \oplus \uvec{0},y)} \right) 
    + (1-\alpha) \left( \phi^B_{(x,\hat{x})} + \phi^A_{(\hat{x},y)}
    \right)\ .
  \end{array}
  $$
  Here, as in the first case, $\Psi_{(x,x \oplus \uvec{0})}$ is
  the flow sending one unit along the arc $(x,x \oplus \uvec{0}) \in
  A$ and nothing along any other arc.
                                
  It is easy to see that this is a valid MCF (i.e. for each pair
  $(x,y)$ the flow $\phi^n_{(x,y)}$ really sends one unit of flow).
  Thus let us check $\phi^n_{\max}$. Firstly, in order to estimate the
  flow on the arcs inside $G[X_A]$ and $G[X_B]$, we determine the
  multiplier by which each flow $\phi^A_{(s,t)}$ respectively
  $\phi^B_{(s,t)}$ appears in the definition of $\phi^n$.  By
  symmetry, it suffices to do this for all pairs $s,t \in X_A$.
                                
  Each pair $s,t \in X'_A$ is used once with multiplier one (for
  $(x,y) = (s,t)$) and once with multiplier $\alpha$ (for $(x,y) = (s
  \oplus \uvec{0},t)$). Thus, each $\phi^A_{(s,t)}$ appears with
  multiplier $(1+\alpha)$ for $s,t\not=\hat{x}$.
  
  Each pair $s = \hat{x}$ and $t \in X'_A$ is used once with
  multiplier one (for $(x,y) = (s,t)$) and, for each of the $(n-1)/2$
  pairs $x \in X_B$ and $y = t$, with multiplier $(1-\alpha)$.
  Each pair $s \in X'_A$ and $t = \hat{x}$ is used once with
  multiplier one (for $(x,y) = (s,t)$) and, for each of the $(n-1)/2$
  pairs $x = s$ and $y \in X_B$, with multiplier $(1-\alpha)$.
                                
  Thus, due to
  $$ 
  (1-\alpha)\frac{n-1}{2} = \frac{n+1 - (n-1)}{n+1} \frac{n-1}{2} =
  \frac{n-1}{n+1} = \alpha\ ,
  $$
  each $\phi^A_{(s,t)}$ with $s = \hat{x}$ or $t = \hat{x}$ appears
  with multiplier $(1+\alpha)$.
  
  Secondly, we estimate the flow along the arcs $(x,x \oplus
  \uvec{0})$ with $x \neq \hat{x}$. By symmetry we restrict our
  attention to the case $x \in X'_A$ and $y \in X'_B$, and we find
  that each arc $(x,x \oplus \uvec{0})$ is used $(n-1)/2$ times with
  flow-value $\alpha$.

  Altogether, this yields
  $$
  \phi^n_{\max} \le \max\left\{ (1+\alpha) \cdot \phi^{n'}_{\max},
    \frac{n-1}{2} \cdot \alpha \right\} \le \frac n2 \ ,
  $$
  which concludes the inductive step, and thus, the proof.
\end{proof}


\section{Towards a Lower-Bound-Theorem for 0/1-Polytopes}

In the following, we will exploit the following construction (using
revlex-initial 0/1-polytopes) several times.

\begin{proposition} \label{prop:pyr}
  For $d,n \in \N$ with $d+1 \le n \le 2^d$ there exists $\tilde{d}
  \in \N$ such that for $\tilde{n} \defi n-(d-\tilde{d})$ the
  following inequalities hold.
  \begin{gather}
    0 \le \tilde{d} \le d 
    \label{eq:pyrtriv} \\
    2^{\tilde{d}-1}\ <\ \tilde{n}\ \le 2^{\tilde{d}}
    \label{eq:lbprop} \\
    \tilde{d}\le 1+\log_2 n \label{eq:log}
  \end{gather}

  Furthermore $\knapP{\tilde{n}}$ is a $\tilde{d}$-dimensional
   revlex-initial 0/1-polytope with $\tilde{n}$ vertices.
\end{proposition}
\begin{proof}
  To see that such a $\tilde{d}$ and $\tilde{n}$ exist, observe that
  with $\tilde{n}(k)\defi n-(d-k)$ we have $\tilde{n}(k) > 2^{k-1}$
  for $k=0$ and $\tilde{n}(k)\le 2^k$ for $k=d$; note that for these
  estimates we need $d+1\le n\le 2^d$. Then, we have that
  $$
  \tilde{d}\ :=\ \min\setdef{k\in\N}{\tilde{n}(k)\le 2^k}
  $$
  satisfies~(\ref{eq:pyrtriv}).
  
  By definition, we have $\tilde{n}(\tilde{d})\le 2^{\tilde{d}}$. If
  $\tilde{d}=0$, then (as stated above) also $\tilde{n}(\tilde{d}) >
  2^{\tilde{d}-1}$ is true, and otherwise, from the minimality
  of~$\tilde{d}$ we conclude $\tilde{n}(\tilde{d}-1)>2^{\tilde{d}-1}$,
  which, of course, implies $\tilde{n}(\tilde{d}) > 2^{\tilde{d}-1}$.
  Hence, $\tilde{d}$ also satisfies~(\ref{eq:lbprop}).

  Finally, (\ref{eq:log}) trivially follows from~(\ref{eq:lbprop}).

  Thus, with $\tilde{n}:=\tilde{n}(\tilde{d})$, by \eqref{eq:lbprop}
  and Proposition \ref{prop:dim} the revlex-initial 0/1-polytope
  $\knapP{\tilde{n}}$ has dimension $\tilde{d}$.
\end{proof}

\begin{definition}
  For arbitrary $d,n \in \N$ with $d+1 \le n \le 2^d$ and 
  $\tilde{d},\tilde{n} \in \N$ as in Proposition \ref{prop:pyr}
  we define $P(d,n)$ to be the $d$-dimensional 0/1-polytope with
  $n$ vertices obtained by building 
  the $(d-\tilde{d})$-fold pyramid over $\knapP{\tilde{n}}$.

  We denote the parameters $\tilde{d}$ and $\tilde{n}$ by
  $\tilde{d}(d,n)$ and $\tilde{n}(d,n)$.
\end{definition}

Note that Proposition \ref{prop:pyr} guarantees that this construction
always works as claimed in the definition of $P(d,n)$.

\subsection {An Upper Bound on the Minimal Number of Facets}

\begin{definition}
  For $d,n\in\N$ with $d+1 \le n \le 2^d$, denote by $\lb{d}{n}$ the
  minimal number of facets of a $d$-dimensional 0/1-polytope with~$n$
  vertices. 
\end{definition}

Note that a $k$-dimensional 0/1-polytope in $\R^d$ (with $k<d$) can
isometrically be projected to a $k$-dimensional 0/1-polytope
in~$\R^k$. Thus, the definition is independent of the ambient spaces
of the polytopes.

\begin{proposition}
\label{prop:lb}
For every $d+1\le n\le 2^d$ we have $\lb{d}{n}\le d+2\log_2 n$.
\end{proposition}

\begin{proof}
  
  By Theorem~\ref{thm:facets}(2), the revlex-initial 0/1-polytope
  $\knapP{\tilde{n}(d,n)}$ has at most $3\tilde{d}-2$ facets. Thus
  $P(d,n)$ has at most $3\tilde{d}-2+n-\tilde{n}=2\tilde{d}+d-2$
  facets.
  The claim of the proposition follows by~\eqref{eq:log}.
\end{proof}

The proposition immediately implies the following results.
\begin{theorem}
\label{thm:lbfac}
\mbox{}
\begin{enumerate}
\item For every $d+1\le n\le 2^d$ we have $\lb{d}{n}\le 3d$.
\item For $d+1\le n(d)\le 2^{\littleO{d}}$ we have $\lb{d}{n(d)}=
  d+\littleO{d}$.
\item For $1<\alpha <2$ and $n(d)\defi\lfloor \alpha^d\rfloor$ we have
  $\lb{d}{n(d)}\le (1+2\log_2\alpha)d+\littleO{d}$.
\end{enumerate}
\end{theorem}

The upper bounds on $\lb{d}{n}$ provided by the polytopes $P(d,n)$ in
Proposition~\ref{prop:lb} are not sharp, at least not for all
parameters~$d$ and~$n$. This follows, for instance, from the examples
of Cartesian products of $r$ 0/1-simplices of dimension
$d_1$,\dots,$d_r$ (which are precisely the simple 0/1-polytopes, see
Kaibel and Wolff~\cite{KW00}). Such a product is a 0/1-polytope of
dimension $d=\sum d_i$ with $\prod(d_i+1)$ vertices and $d+r$ facets.
In particular for $n = (\lfloor\tfrac{d}{2}\rfloor+1)
(\lceil\tfrac{d}{2}\rceil+1)$, this yields
$$
g(d,n)\ =\ d+2\ ,
$$
while the polytopes $P(d,n)$ have $d+\Omega(\log_2 d)$ facets.

\begin{figure}[h]
  \centering
  \includegraphics[width=.45\textwidth]{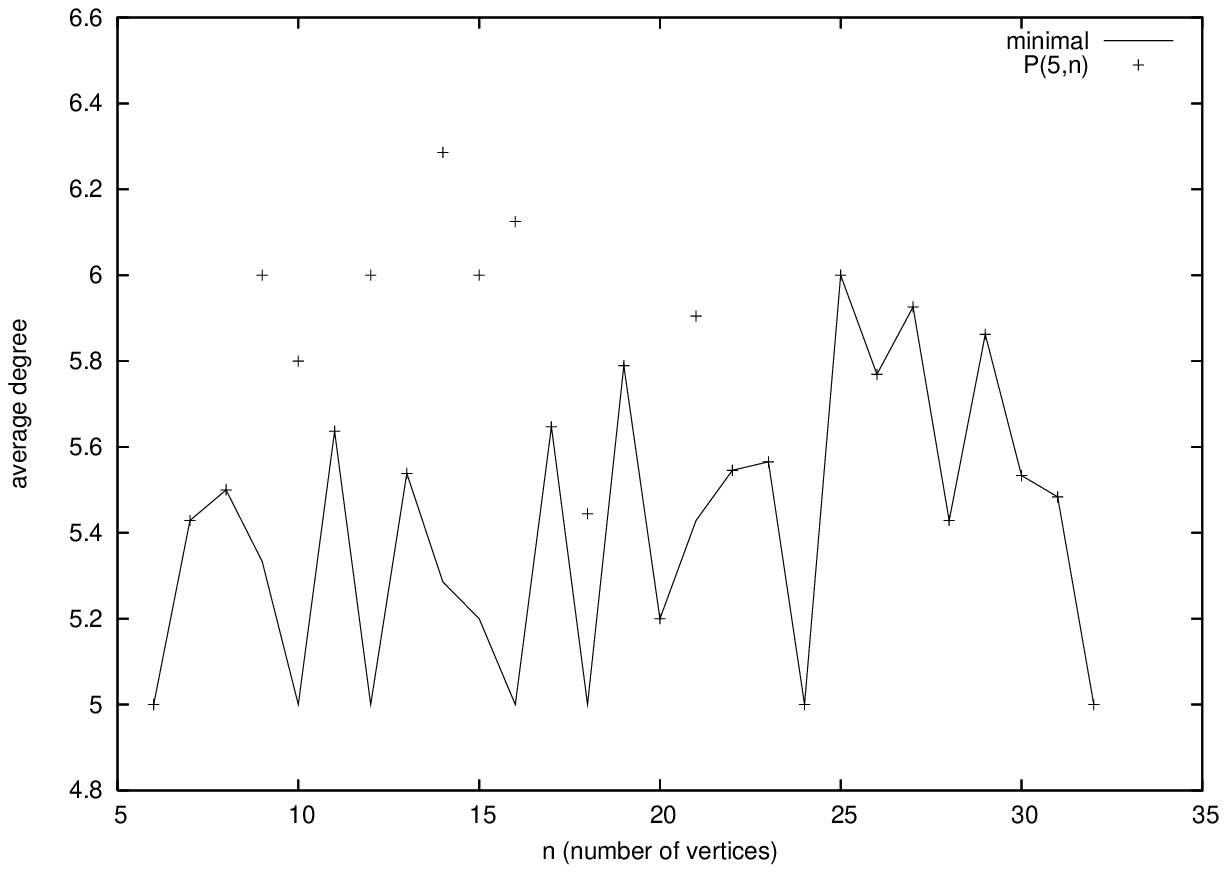}
  \hfill
  \includegraphics[width=.45\textwidth]{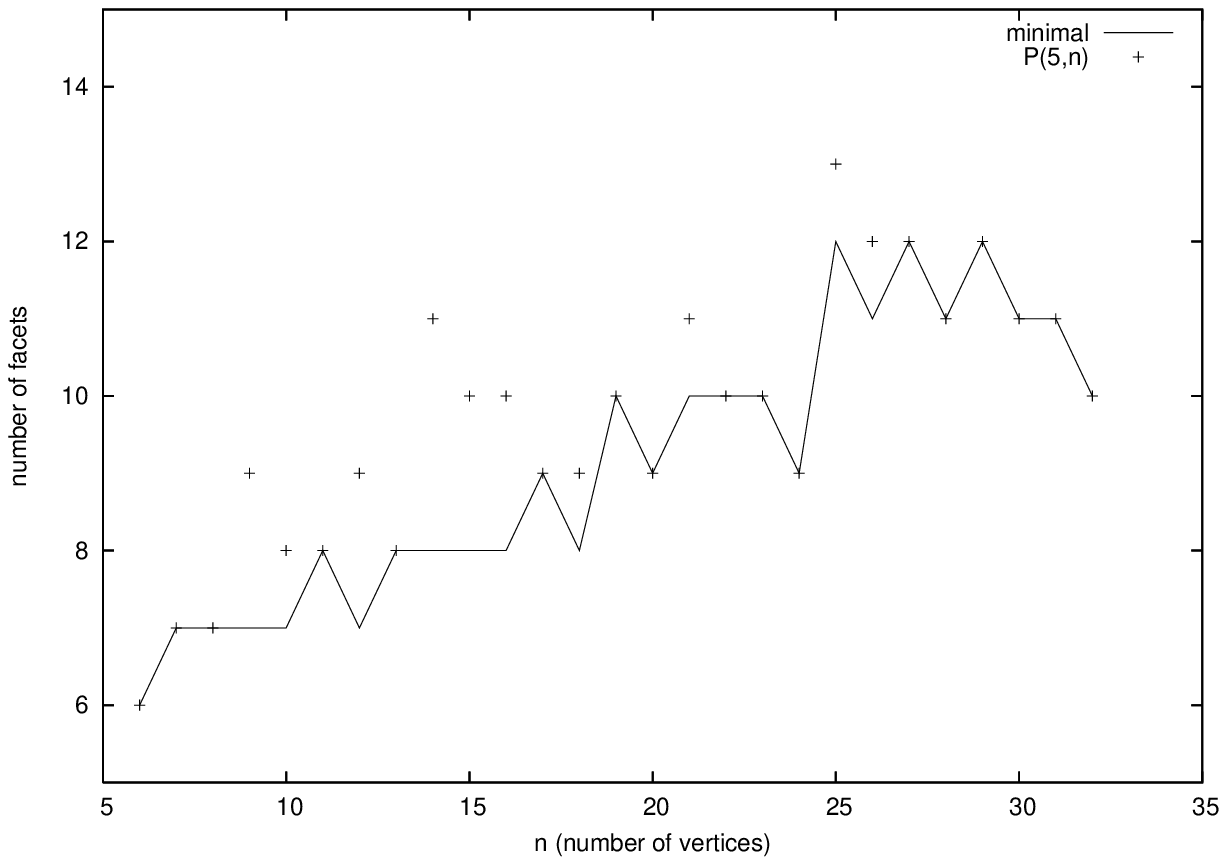}
  \caption{
    Comparison of the lower bounds on $\lb{5}{n}$ and $\lbe{5}{n}$
    obtained from the polytopes $P(5,n)$ in the proofs of
    Propositions~\ref{prop:pyr} and~\ref{prop:lb} with the true values
    of $\lb{5}{n}$ and $\lbe{5}{n}$ obtained from Aichholzer's
    enumeration~\cite{Aic00}.
  }
  \label{fig:aich}
\end{figure}

\begin{figure}[ht]
  \centering
  \includegraphics[width=.6\textwidth]{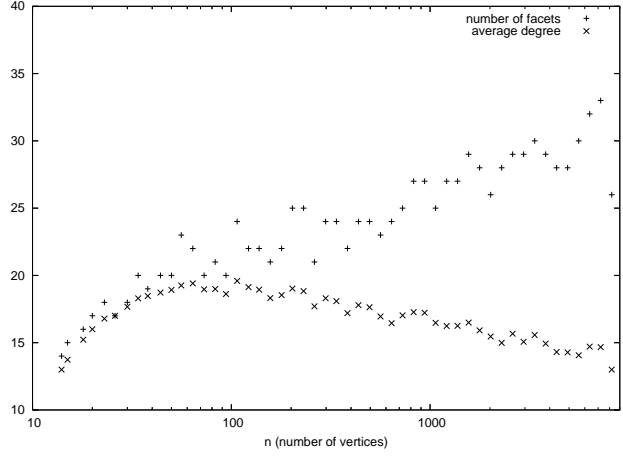}
  \caption{
    Numbers of facets and average degrees of the polytopes $P(13,n)$
    providing the upper bounds on $\lb{13}{n}$ and $\lbe{13}{n}$.  }
  \label{fig:p13}
\end{figure}

The right part of Figure~\ref{fig:aich} shows that for $d=5$ the
polytopes $P(5,n)$ achieve the respective minimum number of facets in
all but 10 cases (i.e., in 17 out of 27 cases).  Figure~\ref{fig:p13}
depicts the numbers of facets (and the average degrees) of the
polytopes $P(13,n)$.

For sub-exponential numbers of vertices, Part~(2) of
Theorem~\ref{thm:lbfac} shows that the minimum number of facets is
asymptotically as small as the number of facets of any $d$-dimensional
polytope can be (up to an additive $\littleO{1}$-term). The range of
sub-exponential vertex numbers is particularly interesting for two
reasons: Firstly, many 0/1-polytopes that are relevant in
combinatorial optimization have sub-exponentially many vertices (e.g.,
cut polytopes of complete graphs and traveling salesman polytopes).
Secondly, the papers by B\'ar\'any and P\'or~\cite{BP01} and
Gatzouras, Giannopoulos, and Markoulakis~\cite{GGM04} show that within
sub-exponential ranges of vertex numbers a random 0/1-polytope has very
many facets. In fact, it may well be that the maximum numbers of
facets of 0/1-polytopes is (roughly) attained by these polytopes.

The examples of products of simplices (i.e., simple 0/1-polytopes)
seem to indicate that it might be hopeless to derive an explicit
formula for $\lb{d}{n}$, i.e., a sharp lower bound theorem for the
facet numbers of 0/1-polytopes. Nevertheless, the question for the
(asymptotic) best upper bound on $\lb{d}{n}$ that does only depend
on~$d$ (and not on~$n$) might be within reach. In particular, we do
not know whether there is some constant $\alpha<3$ such that
$\lb{d}{n}\le \alpha d+\littleO{d}$ holds for all~$d$ and~$n$. This
might even be true for $\alpha=2$.

\subsection{An Upper Bound on the Minimal Number of Edges}

\begin{definition}
  For $d,n\in\N$ with $d+1 \le n \le 2^d$, denote by $\lbe{d}{n}$ the
  minimal average degree among all graphs of $d$-dimensional
  0/1-polytopes with~$n$ vertices.
\end{definition}

Revlex-initial 0/1-polytopes and the pyramidal construction yield the
following bound of the minimum average degrees.

\begin{theorem} \label{thm:edgenumbers}
  For $d+1\le n\le 2^d$, we have $\lbe{d}{n}\le d+4$.
\end{theorem}

\begin{proof} 
  Set $\tilde{d} \defi \tilde{d}(d,n)$ and $\tilde{n} \defi
  \tilde{n}(d,n)$.  By Theorem~\ref{thm:edgenumcompr}, the
  revlex-initial 0/1-polytope $\knapP{\tilde{n}}$ has at most
  $(\tilde{d}+4)\tilde{n}$ edges.  Thus, $P(d,n)$ (the
  $(d-\tilde{d})$-fold pyramid over $\knapP{\tilde{n}}$) has at most
  $$
  (\tilde{d}+4)\tilde{n}+(d-\tilde{d})n\ \le\ (d+4)n
  $$
  edges.
\end{proof}

The left part of Figure~\ref{fig:aich} shows that for $d=5$ the
polytopes $P(5,n)$ achieve the respective minimum average degree in
all but 8 cases (i.e., in 19 out of 27 cases).

\begin{figure}[ht]
  \centering
  \includegraphics[width=.6\textwidth]{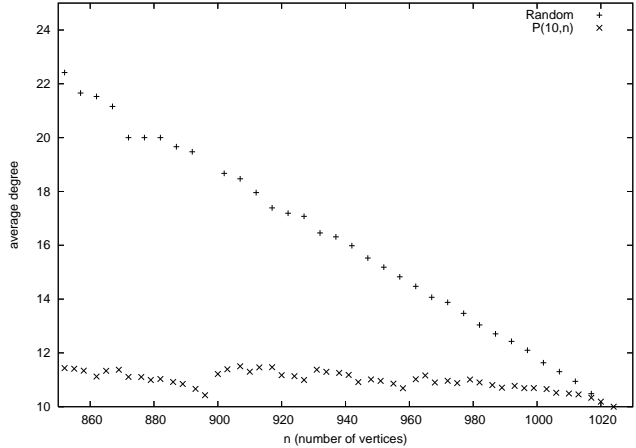}
  \caption{
    Average degrees of the polytopes $P(10,n)$ and uniformly random
    10-dimensional 0/1-polytopes (by sampling).  }
  \label{fig:avg}
\end{figure}

\begin{figure}[ht]
  \centering
  \includegraphics[width=.6\textwidth]{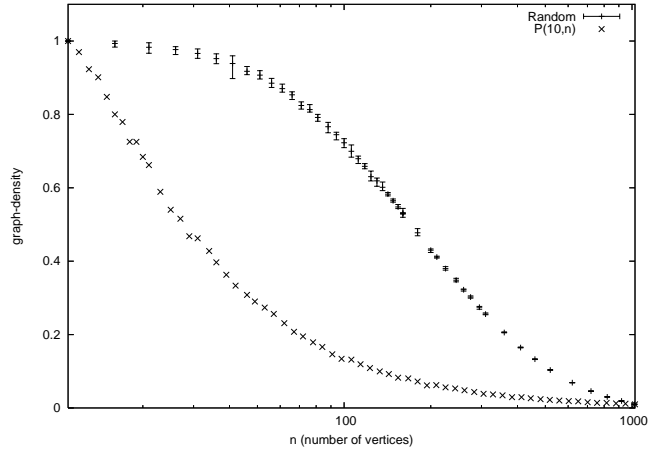}
  \caption{The graph densities of the polytopes $P(10,n)$ used in the
    proof of Theorem~\ref{thm:edgenumbers} versus the graph densities
    of respective random 0/1-polytopes (by sampling).} 
  \label{fig:dens}
\end{figure}


Finally, the polytopes $P(d,n)$ yield examples of 0/1-polytopes with
remarkably sparse graphs, satisfying, nevertheless, the
Mihail-Vazirani conjecture.

\begin{theorem}
\label{thm:expall}
  For every $d+1\le n\le 2^d$, there is a $d$-dimensional 0/1-polytope
  with $n$ vertices, at most $(d+4)n$ edges, and edge expansion at
  least one.
\end{theorem}
\begin{proof}
  By Theorem~\ref{thm:edgenumbers} the polytope $P(d,n)$ has at most
  $(d+4)n$ edges. Since $P(d,n)$ is a $k$-fold pyramid over the
  revlex-initial 0/1-polytope $\knapP{\tilde{n}(d,n)}$ the
  multi-commodity flow constructed in the proof of
  Theorem~\ref{thm:edge-expansion} can be easily extended to a
  multi-commodity flow of $P(d,n)$ sending one unit of flow from every
  vertex to every other vertex.
\end{proof}

\section{Concluding Remarks}

The contributions of this paper concern three topics: (1)
Investigations of a 'natural' class of 0/1-polytopes, (2)
lower bound theorem(s) for 0/1-polytopes, and (3) support of the
Mihail-Vazirani conjecture on the edge expansion of the graphs of
0/1-polytopes.

With respect to the first topic, one may be interested also in
studying the convex hulls of sets of 0/1-vectors that are only
\emph{gradually} revlex-initial (the 0/1-polytopes corresponding to
compressed set systems), i.e., convex hulls of sets~$X$ of 0/1-vectors
which, with every $x\in X$, contain all 0/1-vectors~$y$ which have the
same number of ones as~$x$ and are revlex-smaller than~$x$. Due to the
important role played by the \emph{monotone} ones among them (more
precisely: by the corresponding set systems) in the theory of
simplicial complexes, it might be that these objects bear some
connections between 0/1-polytopes and combinatorial topology.  This
would be quite interesting.

It seems that precise lower bound theorems on the number of facets
(edges, or even other-dimensional faces) are hard to obtain.
Nevertheless, with respect to topic (2) some questions remain open
that may be tractable, e.g., the question whether there is some
$\alpha < 3$ (maybe $\alpha=2$?) with $\lb{d}{n}\le \alpha
d+\littleO{d}$.

Perhaps the most interesting and promising line to follow up this
research concerns topic (3). Extending our techniques for construction
of the multi-commodity flows showing that revlex-initial 0/1-polytopes
(as special knapsack-polytopes) have edge expansion at least one to
all knapsack polytopes (or even to all monotone polytopes) would be a
big support for the Mihail-Vazirani conjecture (which itself is of
great importance in the theory of random generation and approximate
counting, as mentioned in Section~\ref{subsec:edgeexp}). It follows
from work of Morris and Sinclair~\cite{SM99} that the edge-expansion
of the graphs of $d$-dimensional 0/1-knapsack polytopes is bounded
from below by a polynomial in $1/d$. Their proof in fact shows that
this is true even for the subgraph that is formed by those edges which
are also edges of the cube. Since our flows extensively use non-cube
edges, the techniques used in the proof of
Theorem~\ref{thm:edge-expansion} seem to have good potential to
improve the current lower bound, maybe even to 'one' as conjectured by
Mihail and Vazirani.

\bibliographystyle{amsplain}

\begin{thebibliography}{10}

\bibitem{Aic00}
Oswin Aichholzer, \emph{Extremal properties of {$0/1$}-polytopes of dimension
  5}, Polytopes---combinatorics and computation (Oberwolfach, 1997), DMV Sem.,
  vol.~29, Birkh\"auser, Basel, 2000, pp.~111--130.

\bibitem{BP01}
Imre B{\'a}r{\'a}ny and Attila P{\'o}r, \emph{On {$0$}-{$1$} polytopes with
  many facets}, Adv. Math. \textbf{161} (2001), no.~2, 209--228.

\bibitem{FM92}
T.~Feder and M.~Mihail, \emph{Balanced matroids}, Proceedings of the 24th
  Annual ACM ``Symposium on the theory of Computing'' (STOC) (Victoria, British
  Columbia), ACM Press, New York, 1992, pp.~26--38.

\bibitem{FKR00}
Tam{\'a}s Fleiner, Volker Kaibel, and G{\"u}nter Rote, \emph{Upper bounds on
  the maximal number of facets of 0/1-polytopes}, European J. Combin.
  \textbf{21} (2000), no.~1, 121--130, Combinatorics of polytopes.

\bibitem{GGM04}
D.~Gatzouras, A.~Giannopoulos, and N.~Markoulakis, \emph{Lower bound for the
  maximal number of facets of a 0/1 polytope}, Tech. report, University of
  Athens, 2004, To appear in \emph{Discrete Comp. Geom.}

\bibitem{JS97}
M.~Jerrum and A.~Sinclair, \emph{The {M}arkov {C}hain {M}onte {C}arlo method:
  An approach to approximate counting and integration}, Approximation
  Algorithms (D.~Hochbaum, ed.), PWS Publishing Company, Boston, 1997,
  pp.~482--520.

\bibitem{Kai04}
Volker Kaibel, \emph{On the expansion of graphs of 0/1-polytopes}, The Sharpest
  Cut: The Impact of Manfred Padberg and His Work (Martin Gr\"otschel, ed.),
  MPS-SIAM Series on Optimization, vol.~4, SIAM, 2004, pp.~199--216.

\bibitem{KW00}
Volker Kaibel and Martin Wolff, \emph{Simple 0/1-polytopes}, European J.
  Combin. \textbf{21} (2000), no.~1, 139--144.

\bibitem{Mih92}
M.~Mihail, \emph{On the expansion of combinatorial polytopes}, Proceedings of
  the 17th International Symposium on ``Mathematical Foundations of Computer
  Science'' (I.~M. Havel and V.~Koubek, eds.), Lecture Notes in Computer
  Science, vol. 629, Springer-Verlag, 1992, pp.~37--49.

\bibitem{SM99}
B.~Morris and A.~Sinclair, \emph{Random walks on truncated cubes and sampling
  0-1 knapsack problem}, Proceedings of the 40th IEEE Symp. on Foundations of
  Computer Science (New York), 1999, pp.~230--240.

\bibitem{Sch86}
Alexander Schrijver, \emph{{Theory of linear and integer programming.}}, {John
  Wiley \& Sons}, 1986.

\bibitem{Sin93}
A.~Sinclair, \emph{{Algorithms for random generation and counting: a Markov
  Chain approach.}}, Progress in Theoretical Computer Science, {Birkh\"auser},
  Boston, 1993.

\bibitem{Sta80}
Richard~P. Stanley, \emph{Decompositions of rational convex polytopes}, Ann.
  Discrete Math. \textbf{6} (1980), 333--342.

\bibitem{Zie98}
G{\"u}nter~M. Ziegler, \emph{Lectures on polytopes}, Graduate Texts in
  Mathematics, vol. 152, Springer-Verlag, New York, 1995, Revised edition:
  1998.

\bibitem{Zie00}
\bysame, \emph{Lectures on {$0/1$}-polytopes}, Polytopes---combinatorics and
  computation (Oberwolfach, 1997), DMV Sem., vol.~29, Birkh\"auser, Basel,
  2000, pp.~1--41.

\end{thebibliography}

\providecommand{\bysame}{\leavevmode\hbox to3em{\hrulefill}\thinspace}
\providecommand{\MR}{\relax\ifhmode\unskip\space\fi MR }
\providecommand{\MRhref}[2]{%
  \href{http://www.ams.org/mathscinet-getitem?mr=#1}{#2}
}
\providecommand{\href}[2]{#2}

\end{document}